\documentclass[12pt]{article}
\usepackage{a4,verbatim,amsmath,amssymb}
\usepackage[isolatin]{inputenc}
\usepackage{epsfig}
\input math1.def
\usepackage{amsthm}

\newtheorem{theor}{Theorem}[section]

\newtheorem{lemma}[theor]{Lemma}

\newtheorem{prop}[theor]{Proposition}

\numberwithin{equation}{section}
\def\proof{\goodbreak\noindent{\sc Proof. }\nobreak}

\def\endproof{\par\nobreak\hbox to \hsize{\hfil\vrule width 5pt height
5pt}\goodbreak\vskip 3pt}

\begin{document}
\title{ Conditional Square Functions and Dyadic Perturbations of the Sine-Cosine decomposition
for Hardy Martingales}
     
\author{Paul F. X. M\"uller\thanks{Supported
by the Austrian Science foundation (FWF) Pr.Nr. FWFP28352-N32.  }}
\date{October  $16$,  2016}
\maketitle
\begin{abstract} 
We prove that the $\cP -$norm estimates between a Hardy martingale and its cosine part
are stable under dyadic perturbations.
\end{abstract}
\tableofcontents
\section{Introduction}\label{intro}
Hardy martingales developed alongside Banach spaces of analytic functions and
played an important role in establishing  their isomorphic invariants. 
For instance those martingales were employed  in the construction of subspaces in $L^1/H^1$ 
isomorphic to $L^1.$  An integrable Hardy martingale $F=(F_k)$ satisfies 
the $L^1$ estimate 
$$ \| \sup_k |F_k| \| _1 \le e   \sup_k \|F_k  \| _1 , $$ 
and it may be decomposed into the sum of Hardy martingales as $ F = G+B$ such that 
$$ \| (\sum \bE_{k-1} |\Delta _k G|^2)^{1/2} \|_1 + \sum \|\Delta B_k\|_1 \le C\|F\|_1. $$
See Garling, Bourgain, Mueller. Equally peculiar  for Hardy martingales are the 
are the transform estimates
$$
 \| (\sum \bE_{k-1} |\Delta _k G|^2)^{1/2} \|_1 \le C  
\| (\sum \bE_{k-1} |\Im w_{k-1}\Delta _k G|^2)^{1/2} \|_1 ,$$
for  every  adapted sequence $(w_k) $ satisfying $ |w_k| \ge 1/C . $ 
A proof of Bourgain's theorem that $L^1$ embeds into   $L^1/H^1$ may be obtained 
in the following way: 
\begin{enumerate} 
\item Use as starting point the  estimates of the Garnett Jones Theorem.
\item Prove stability under dyadic perturbation for the Davis and Garsia Inequalities.
\item Prove stability under dyadic perturbation of  the martingale transform estimates. 
\end{enumerate}
We determined the extent to which DGI are stable under dyadic perturbation, and we showed how 
the above strategy actually gives an isomorphism from $L^1$ into a subspace of $L^1/H^1. $ 
In the present paper we turn to the martingale transform estimates and verify that they are 
indeed stable under dyadic perturbations. 
\section{Preliminaries}

\paragraph{Martingales and Transforms on  $\bT^\bN$. }
Let  $\bT=  \{ e^{i\theta} : \theta \in
[0, 2\pi [ \} $ be  the torus 
equipped with the  normalized angular measure.
 Let $\bT^\bN $
be its  countable product
equipped with the  product Haar measure $\bP .$ We let 
$\bE $ denote expectation  with respect to   $\bP .$ 

Fix $k \in \bN $, the 
cylinder sets 
$ \{(A_1, \dots, A_k , \bT^\bN )\},$
where $A_i,\, i \le k $ are measurable subsets 
of $\bT$,  form the $\s-$algebra  $\cF_k $.
Thus we obtain a filtered probability space $(\bT^\bN, (\cF_k) , \bP) $.
We let $\bE_{k}$ denote the conditional expectation with
respect to the $\s-$algebra  $\cF_k .$
Let  $G = (G_k) $ be  an $L^1(\bT^\bN)-$bounded  martingale.
Conditioned on  $\cF_{k-1}$  the martingale difference $\Delta G_k =
G_k- G_{k-1}$  defines an element in  $L_0^1(\bT) ,$ the Lebesgue space of integrable,
 functions with vanishing mean.
We define the  previsible norm  as 
\begin{equation}\label{15augmt2}
\| G\|_\cP  =\| ( \sum_{k = 1 }^\infty \bE_{k-1} |\Delta G_k |^2 )^{1/2} \|_{L^1},
\end{equation}
and  refer to $( \sum_{k = 1 }^\infty \bE_{k-1} |\Delta G_k |^2 )^{1/2}$
as the conditional square function of $G.$  

For    any bounded and 
adapted sequence  $W = ( w_k ) $  we define   
 the martingale transform  operator $T_W $ by 
\begin{equation}\label{18-10-19-4}
T_W (G )  = \Im \left [ \sum w_{k-1} \Delta_k G  \right ] . 
\end{equation}

  Garsia \cite{sia}  is our reference to  martingale inequalities.

\paragraph{Sine-Cosine decomposition.}

Let 
 $G = (G_k)
$ 
be a  martingale on $\bT^\bN$ with respect to the canonical product filtration   $(\cF_k) $. 
Let   $U = (U_k)$  
be the martingale  defined by averaging
\begin{equation}\label{18-10-16-2}
 U _k (x,y)   = 
\frac12 \left[  G_k(x,y)+   G   _k(x,\overline{y})\right], 
\end{equation}
%
where $x \in \bT^{k-1} , \,  y \in \bT .   $ 
The martingale   $U $ 
is called the cosine part of $G$. 
Putting $V_k = G_k - U_k $ we obtain the corresponding 
sine-martingale $ V = ( V_k) $, and the  sine-cosine decomposition of $G$ defined by 
$$ G = U +V . $$  
By construction we have 
$\Delta V_k ( x , y ) = -\Delta V_k ( x , \overline{y} ), $ 
and $  U_k ( x , y ) = U_k ( x , \overline{y} ), $ for any $k\in \bN .$ 
\paragraph{The Hilbert transform.}
The Hilbert transform on $L^2 ( \bT )$ is defined as 
Fourier multiplier by 
$$ H(  e^{in\theta}) = -i {\rm sign} (n)  e^{in\theta} . $$

Let $1 \le p \le \infty. $ The Hardy space    $H^p_0 (\bT ) \sb  L^p_0 (\bT ) $   
consist of those $p-$integrable
functions of vanishing mean,  for which the harmonic extension to the
unit disk is analytic. See  \cite{jg81}.
For  $ h \in H^2_0( \bT )$ and let $ y = \Im h . $
The Hilbert transform recovers $h $ from its imaginary part $y$ ,
we have $ h = -Hy +iy . $ and 
$ \| h \| _2 = \sqrt{2} \| y \|_2 . $ For 
$ w \in \bC ,$ $ | w| = 1 $ we have therefore
$$   \| h \| _2  =    \sqrt{2}  \| y \| _2= \sqrt{2}\| \Im (w\cdot h)\|_2 . $$

\section{Martingale estimates}
\paragraph{Hardy martingales. }

An  $L^1(\bT^\bN ) $ bounded  $(\cF_k)$ martingale $G = (G_k) $
is called a Hardy martingale if conditioned on $\cF_{k-1}$ 
the martingale difference  
$ \Delta G _k $ 
defines an element in $H^1_0 (\bT ).  $  See \cite{g1}, \cite{gar2}.
\cite{pfxm12, deco-2014, deco-2016}

Since  the Hilbert transform,  applied to  functions with vanishing mean,    preseves   the   $ L^2 $ norm, we have      
$\bE_{k-1} | \Delta U _k |^2 
=  \bE_{k-1} |\Im w_{k-1} \Delta G_k|^2 , $
 for each adapted sequence  $W = ( w_k ) $  with  $|w_k | =1 ,$ and consequently, 
\begin{equation}\label{24o1210}
\| ( \sum   \bE_{k-1} | \Delta U_k |^2  )^{1/2} \|_1 = 
 \| ( \sum   \bE_{k-1} |\Im w_{k-1} \Delta G_k|^2   )^{1/2} \|_1 .
 \end{equation}
We  restate  \eqref{24o1210} as  
 $  
  \| U \| _{\cP }  =  \| T_ W ( G ) \|_\cP ,$
where 
$
T_W (G  )  = \Im \left [ \sum w_{k-1} \Delta_k (G  )\right ] . 
$
In this paper we show that the lower $\cP$ norm estimate $ 
  \| U \| _{\cP }  \le   \| T_ W ( G ) \|_\cP ,$ is  stable under dyadic perturbation. 

\paragraph{Dyadic martingales. }
The dyadic sigma-algebra on $ \bT^\bN $ is  defined with 
 Rademacher functions. For $ x = ( x _k ) \in  \bT^\bN $
define $\cos_k ( x ) =  \Re x_k  $ and 
$$ 
 \s_ k ( x ) = {\rm sign} (\cos_k ( x )).
$$
We let  $\cD$ be the  sigma- algebra 
 generated by 
$ \{\s_k , k \in \bN \}  $ and call it the dyadic sigma-algebra
on $ \bT^\bN .$ 
Let  $ G \in L^1 ( \bT^\bN ) $ 
with  sine cosine decomposition $G = U + V $, 
then  $\bE ( U_k  | \cD )  =\bE ( G _k | \cD ) $
for $k\in \bN , $ and hence 
$$ U -  \bE ( U  | \cD ) + V = G -  \bE ( G  | \cD ). $$ 
Our principle result 
asserts    
 stability  for \eqref{24o1210} under  dyadic perturbations as follows:  

\begin{theor}\label{11aug60b}
Let  $G = (G_k)_{k = 1 }^n $ be a  martingale and let  $U = (U_k)_{k = 1 }^n $
be its cosine martngale given by \eqref{18-10-16-2}. 
Then, for    any 
adapted sequence  $W = ( w_k ) $ satisfying   $|w_k | =1  ,$ we have  
\begin{equation}\label{11aug95aa}
\|U - \bE ( U |\cD ) \|_{\cP} \le 
C \| T_ W ( G - \bE ( G | \cD ))  \|_\cP^{1/2} \|G\|_{\cP}^{1/2} ,\end{equation}
where $T_W $ is the martingale transform  operator defined by \eqref{18-10-19-4}.
\end{theor}

Define   $\s \in L^2 (\bT )$ by  
$ \s (\zeta) = {\rm sign} \Re \zeta . $ Note that $\s (\zeta) = \s(\overline{\zeta}) , $
for all $\zeta \in \bT . $  
For $ f , g \in L^2(\bT) $ we put $ \la f , g \ra =\int_\bT f \overline {g} dm .$  
\begin{lemma}\label{10sep121}
Let   $h \in H_0^2 ( \bT) ,$ and   
$  u(z) 
= (h(z) +h( \overline{z}  ))/2 $.
Then for $ w,b \in \bC , $ with  $ |w| = 1 , $
$$
 \Im ^2(w\cdot(\la u, \s \ra - b))  + \Re^2 (w\cdot\la u, \s \ra) + 
\int_{\bT} |u - \la u, \s \ra \s |^2 d m 
=\int_{\bT} \Im ^2(w\cdot( h - b \s )) d m
$$
\end{lemma}
\proof 
First   put 
$ w_0 = 1_\bT , $ $ w_1 = \s, $ 
and choose any    orthonormal system $\{w_k: k \ge 2 \}$ in   
$L_G^2 (\bT ) $ so that  $\{w_k: k \ge 0\}$ is an orthonormal
basis for  $L_G^2 (\bT ). $ Then  
$\{w_k, H w_k : k \ge 0\},$
where $H$  the Hilbert transform, is a orthonormal basis in $ L^2 ( \bT )$.
Moreover in the  Hardy space  $ H^2 ( \bT )$  the  analytic system 
$$
\{(w_k +i H w_k) : k \ge 0\}$$
 is an  orthogonal
basis with $\|w_k +i H w_k \|_2 = \sqrt{2} , \, k \ge 1 .$

Fix $ h \in H_0^2 ( \bT) $ and  $ w,b \in \bC , $ with  $ |w| =
1  .$  Clearly by replacing $ h $ by $w h $ and $ b $ by $ w b $ 
it suffices to prove the lemma  with  $w = 1 . $
Since $ \int u = 0 $ we have
that $$ u =   \sum_{n = 1}^\infty c_n w_n .$$ 
We apply the Hilbert transform and  rearrange terms to get  
\begin{equation}\label{10sep122}
 h - b\s= (c_1 - b)\s + ic_1 H\s + \sum_{n = 2}^\infty c_n (w_n +i H
 w_n) . 
\end{equation}
Then, taking imaginary parts gives 
\begin{equation}\label{10sep123} 
\Im (h-b\s) = 
\Im (c_1 - b)\s  + \Re c_1 H \s   + \sum_{n = 2}^\infty  \Im c_n w_n + 
\Re c_n Hw_n . 
\end{equation}
By ortho-gonality the  identity \eqref{10sep123} yields 
\begin{equation}\label{10sep124}
\int_{\bT} \Im ^2( h - b\s ) d m
= \Im ^2(c_1 - b)  + \Re^2 c_1 + \sum_{n = 2}^\infty  |c_n|^2 .
\end{equation}
On the other hand, since  $ \int u = 0 $, $c_1 =\la u, \s \ra, $  and $ w_1 = \s $
we get 
\begin{equation}\label{10sep125}
\int_{\bT} |u - \la u, \s \ra \s |^2 d m
 =  \sum_{n = 2}^\infty |c_n|^2 .
\end{equation} 
Comparing the equations \eqref{10sep124}  and \eqref{10sep125}
completes the proof. 
\endproof


We use below  some arithmetic, that we isolate first.
\begin{lemma}\label{10sep126}
Let  $ \mu , b \in \bC $ and
 \begin{equation}\label{11aug65d}
 |\mu| +\dfrac{|\mu-b|^2}{|\mu|+|b|} = a. 
\end{equation}
Then for any  $w \in \bT ,$ 

\begin{equation}\label{11aug60e}
( a - |b|)^2 \le 4( \Im^2 (w\cdot(\mu-b))
+ \Re^2 (w\cdot \mu) )         .
\end{equation}
and 
\begin{equation}\label{11aug60d}
|\mu-b |^2 \le 2(a^2 - |\mu|^2). 
\end{equation}

\end{lemma}
\proof 
By rotation invariance it suffices to prove 
\eqref{11aug60e} for $w = 1 .$ 
Let $ \mu = m_1 + i m_2 $ and $ b = b_1 + i b_2 . $
By definition \eqref{11aug65d}, we have 
$$ 
a - |b| = \dfrac{ |\mu| ^2 - |b|^2 + | \mu - b |^2}{| \mu| + |b |}.
$$ 
Expand and regroup the numerator
\begin{equation}\label{17sep121}
|\mu| ^2 - |b|^2 + | \mu - b |^2 
 = 2m_1 ( m_1  - b_1 ) + 2 m_2 ( m_2 - b_2 ) .
\end{equation}
By the Cauchy Schwarz inequality, the right hand side  \eqref{17sep121}
is bounded by 
$$ 
2 ( m_1^2 + ( m_2 - b_2 )^2 )^{1/2}(  m_2^2 + ( m_1 -b_1)^2)^{1/2} .  
$$
Note that  $m_1 = \Re \mu $   and  $m_2 - b_2 = \Im ( \mu - b ) . $ 
It remains to observe that 
$$
(  m_2^2 + ( m_1 -b_1)^2)^{1/2}  \le |\mu | + |b| .$$
or equivalently 
$$
 m_1^2+ m_2^2  - 2 m _1 b_2 + b_1^2  \le  |\mu |^2 +2   |\mu| |b|+
 |b|^2, $$
which is obviously true. 

Next we turn to verifying \eqref{11aug60d}.  We have
$a^2 -|\mu|^2 = (a+|\mu|)(a -|\mu|)$ hence
\begin{equation}\label{17sep124}
a^2 -|\mu|^2 = 
\left[ 2|\mu| +\dfrac{|\mu-b|^2}{|\mu|+|b|}
\right]\dfrac{|\mu-b|^2}{|\mu|+|b|} .\end{equation}
In view of \eqref{17sep124} we get   \eqref{11aug60d} by showing that  
\begin{equation}\label{17sep123}
 2|\mu|^2 +2|\mu||b| +  |\mu-b|^2 \ge \frac12 ( |\mu|+|b| )^2 .
\end{equation}
The left hand side of \eqref{17sep123} is larger than 
$|\mu|^2+|b|^2$ while the right hand side of   \eqref{17sep123} is
smaller $|\mu|^2+|b|^2 .$ 
\endproof

We  merge the inequalities of Lemma \ref{10sep126} 
with the identity in Lemma \ref{10sep121}.


\begin{prop}\label{25a121} 
 Let $ b \in \bC  $  and $h \in H_0^2 ( \bT) .$
 If  $  u(z) 
= (h(z) +h( \overline{z}  ))/2 $ and 
 $$  | \la u, \s \ra| +\frac{| \la u, \s \ra-b|^2}{| \la u,
  \s \ra| + |b|}   = a , $$
then 
\begin{equation}\label{22aug12b}
\int_{\bT} |u - b \s |^2 dm
\le   8(a^2 -|\la u, \s \ra|^2)+ \int_{\bT} |u - \la u, \s \ra \s |^2
dm.
\end{equation}
and  for all  $ w  \in \bC , $  with $|w| = 1 ,$    
\begin{equation}\label{22aug12}
( a - |b| )^2 + \int_{\bT} |u - \la u, \s \ra \s |^2
dm \le 8 
\int_{\bT}  \Im ^2(w\cdot( h - b \s )) dm .                           
\end{equation}  
\end{prop}

\proof
Put 
\begin{equation}\label{25a123}
 J^2 = \int_{\bT}  \Im ^2(w\cdot( h - b \s )) dm .  
\end{equation}
The proof  exploits the  basic identities for the integral $J^2 $ and 
$\int_{\bT} |u - b \s |^2 dm $ and intertwines them   with  the
arithmetic    \eqref{11aug65d}  -- \eqref{11aug60d}.
\paragraph{Step  1.} Use  the  straight forward identity,
\begin{equation}\label{22aug12a}
\int_{\bT} |u - b \s |^2 dm
= |\la u, \s \ra   -b |^2 + \int_{\bT} |u - \la u, \s \ra \s |^2
dm.
\end{equation}
Apply \eqref{11aug60d},  so that  
$$  |\la u, \s \ra   -b |^2 \le 8(a^2 -|\la u, \s \ra|^2 ) ,$$
hence by \eqref{22aug12a} we get \eqref{22aug12b},
\begin{equation*}
\int_{\bT} |u - b \s |^2 dm
\le   8(a^2 -|\la u, \s \ra|^2)+ \int_{\bT} |u - \la u, \s \ra \s |^2
dm.
\end{equation*}
      
\paragraph{Step  2.} 
The identity of Lemma  \ref{10sep121}  gives 
\begin{equation}\label{22aug12c}
 \Im ^2(w\cdot(\la u, \s \ra - b))  + \Re^2 (w\cdot\la u, \s \ra) + 
\int_{\bT} |u - \la u, \s \ra \s |^2 dm = J^2 .
\end{equation}
Apply \eqref{11aug60e}  with $ \mu = \la u, \s \ra  $
to the left hand side in \eqref{22aug12c}, and
get \eqref{22aug12},
\begin{equation*}
( a - |b| )^2 + \int_{\bT} |u - \la u, \s \ra \s |^2
dm \le 8 J^2 .
\end{equation*}  

\endproof

\subsubsection*{Proof of Theorem \ref{11aug95aa}}
 Let $ \{g_k\} $  be the martingale difference sequence of 
the Hardy martingale $G = (G_k) ,$ and 
let    
 $ \{u_k\} $  be the martingale difference sequence of
the  associated cosine martingale  
$U = (U_k) .$ 
 By convexity  we have 
$$\bE  (  \sum_{k = 1}^\infty |\bE_{k-1}( u_k  \s_k )|^2)^{1/2}  =  
 \bE \bE ( (  \sum_{k = 1}^\infty |\bE_{k-1}( u_k  \s_k )|^2)^{1/2}  | \cD )  \ge \bE  (  \sum_{k = 1}^\infty |\bE (\bE_{k-1}( u_k  \s_k ) | \cD) |^2
)^{1/2}  .
$$
Put  $ b_k  = \bE (\bE_{k-1}( u_k  \s_k ) | {\cD})$
and note that   
$   \bE ( u_k | \cD ) =  
b_k \s_k
.$ 
\paragraph{Step 1.}
Let  $
Y^2 = \sum _{k = 1 }^\infty | \bE_{k-1} (u_k  \s_k) |^2 $ and $Z^2 = \sum_{k = 1 }^\infty  |b_k| ^2  $. 
Then  restating the  above  convexity estimate we have 
%
\begin{equation}\label{28aug70j}
\bE ( Y ) \ge \bE ( Z) .
\end{equation} 
\paragraph{Step 2.} 
Since  
 $
\bE(g_k |{\cD}) = \bE(u_k  | {\cD}) , $
the square of the conditioned square functions of 
$ T_W (   G -  \bE ( G | {\cD})) $  coincides with
\begin{equation}\label{28a12x} 
\sum  \bE_{k-1} |\Im (w_{k-1}\cdot( g_k  - b_k \s_k ))|^2 .
\end{equation}
\paragraph{Step 3.} 
The sequence $ \{ u_k -b_k \s_k\} $  is the  martingale difference sequence of  
$ U - \bE_{\cD}( U)$.
The square of its  conditioned square functions  is hence given by  
\begin{equation}\label{28a1245}
\sum \bE_{k-1}| u_k - b_k \s_k |^2 .
\end{equation}
Following the pattern of \eqref{11aug65d} define 
$$ a_k  = |\bE_{k-1}( u_k \s_k )| +\frac{|\bE_{k-1}( u_k  \s_k )   -b_k|^2}{ 
|\bE_{k-1}( u_k \s_k )| +   |b_k|} 
 , $$
and   
$$v_k =u_k - \bE_{k-1}( u_k  \s_k )\s_k , \quad\quad  r_k ^2 = \bE_{k-1}| v_k  |^2 .$$ 
 By  \eqref{22aug12b} 
\begin{equation}\label{11aug70a}
\begin{aligned}
&\bE_{k-1}| u_k - b_k \s_k |^2 
\le  8 ( a_k^2  +  r_k ^2  - |\bE_{k-1}^2( u_k  \s_k ) |     ).
\end{aligned}
\end{equation}

\paragraph{Step 4.}
With  
$ X^2 = \sum_{k = 1 }^\infty a_k ^2 + r_k^2 ,$
we have the obvious pointwise estimate,  $ X \ge Y $. Taking into account \eqref{11aug70a} 
gives 
 \begin{equation}\label{28aug70g}
 \| U - \bE (U | {\cD}) \|_\cP \le \sqrt{8}  \bE (X^2 -Y^2)^{1/2} \le \sqrt{8}( \bE( X - Y ))^{1/2} ( \bE( X + Y ))^{1/2} . 
\end{equation}
 The  factor $\bE( X + Y )$ 
in \eqref{28aug70g}  admitts an upper bound by 
\begin{equation}\label{28aug70h}
 \bE (X+Y) \le C \|U  \|_{\cP} \le  C \|G  \|_{\cP} .
\end{equation}
\paragraph{Step 5.}
Next we turn to  estimates for  $\bE( X - Y ).$ 
By \eqref{28aug70j}, 
$\bE( X - Y ) 
\le 
 \bE( X - Z ) ,$
and by triangle inequality
$$
X - Z \le (\sum_{k = 1 }^\infty (a_k - |b_k|)^2 + r_k^2)^{1/2} 
.$$
By \eqref{22aug12} 
\begin{equation*}
( a_k - |b_k|)^2 +  r_k^2 
\le 8 \bE_{k-1} |\Im (w_{k-1}\cdot( g_k  - b_k \s_k ))|^2 ,
\end{equation*} and hence 
$$ \bE( X - Z ) 
 \le C \| T_W(G  - \bE(  G | {\cD}))\|_\cP .$$
Invoking 
\eqref{28aug70g}
and \eqref{28aug70h} completes the proof.
\endproof

\bibliographystyle{abbrv}
\bibliography{hardymartingales}
Department of Mathematics\\
J. Kepler Universit\"at Linz\\
A-4040 Linz\\
paul.mueller@jku.at

\end{document}